\documentclass[11pt, reqno]{amsart}
\usepackage{amssymb}
\usepackage{amsmath}
\usepackage{amscd}
\usepackage{amssymb}
\usepackage{eucal}
\usepackage{amsmath}
\usepackage{amscd}
\usepackage[dvips]{color}
\usepackage{multicol}
\usepackage[all]{xy}
\usepackage{graphicx}
\usepackage{color}
\usepackage{colordvi}
\usepackage{xspace}
\usepackage{tikz}

\usepackage{ifpdf}
\ifpdf
 \usepackage[colorlinks, linkcolor=blue, final,backref=page,hyperindex]{hyperref}
\else
 \usepackage[colorlinks, linkcolor=blue, final,backref=page,hyperindex,hypertex]{hyperref}
\fi

\allowdisplaybreaks
\newtheorem{theorem}{Theorem}[section]

\newtheorem{lem}[theorem]{Lemma}
\newtheorem{coro}[theorem]{Corollary}
\newtheorem{prop-def}{Proposition-Definition}[section]

\newtheorem{remark}[theorem]{Remark}

\newenvironment{pf}{\noindent\textit{Proof.}}{\hspace*{1mm} \hfill$\Box$}

\newcommand{\Reg}{\mathrm{reg}}
\newcommand{\Pd}{\mathrm{pd}}
\newcommand{\Max}{\mathrm{max}}

\newcommand{\ra}{\mathrm{\longrightarrow }}
\newcommand{\K}{\mathbb{K}}

\newcommand{\nc}{\newcommand}
\nc{\tred}[1]{\textcolor{green}{#1}}
\nc{\tblue}[1]{\textcolor{blue}{#1}}
\nc{\tgreen}[1]{\textcolor{green}{#1}}
\nc{\tpurple}[1]{\textcolor{purple}{#1}}
\nc{\btred}[1]{\textcolor{green}{\bf #1}}
\nc{\btblue}[1]{\textcolor{blue}{\bf #1}}
\nc{\btgreen}[1]{\textcolor{green}{\bf #1}}
\nc{\btpurple}[1]{\textcolor{purple}{\bf #1}}

\nc{\Jw}[1]{\textcolor{green}{JW: #1}}
\nc{\To}[1]{\textcolor{blue}{#1}}

\begin{document}
\title{Regularity of powers of path ideals of line graphs}


\author {Jiawen Shan}
\address{School of Mathematical Sciences,  Soochow University, 215006 Suzhou, P.R.China}
\email{ShanJiawen826@outlook.com}
\author {Dancheng Lu$^*$}
\address{School of Mathematical Sciences,  Soochow University, 215006 Suzhou, P.R.China}
\email{ludancheng@suda.edu.cn}
\dedicatory{Dedicated to the memory of Professor J\"urgen Herzog }

\footnote{Corresponding author}

\maketitle
\baselineskip=18pt

\begin{abstract} Let $I$ be the $t$-path ideal  of $L_n$, the line graph with $n$-vertices.  It is shown  that $I^s$ has a linear resolution for some $s\geq 1$ (or equivalently for all $s\geq 1$) if and only if $I^s$ has linear quotients for some $s\geq 1$ (or equivalently for all $s\geq 1$) if and only if $t\leq n\leq 2t$. In addition, we present an explicit formula for the regularity of $I^s$ for all $s\geq 1$. It turns out it is linear in $s$ from the very beginning.
\end{abstract}

\section{Introduction}
\setcounter{equation}{0}
\renewcommand{\theequation}
{1.\arabic{equation}}
Let $R=\K[x_1,\ldots, x_n]$ be the polynomial ring in variables $x_1,\ldots, x_n$ over a field $\K$. For a homogeneous ideal $I\subset R$, the Castelnuovo-Mumford regularity (regularity for short) of $I$ is an important algebraic invariant which measures its complexity. The study on the regularity of powers of homogeneous ideals began with a celebrated result, which was independently proved by Cutkosky-Herzog-Trung and Kodiyalam in \cite{CHT99} and \cite{Kod00} respectively. This result states that the regularity of $R/I^s$ is asymptotically a linear function in $s$, i.e., there exist integers $d, e$ and $s_0$ such that $\Reg~R/I^s=ds+e$ for all $s\geq s_0$. Scholars' research in this respect can be roughly divided into two classes. The one is to understand the constant $e$ and the minimum power $s_0$ starting from which $\Reg~R/I^s$ becomes a linear function for some special or general ideals $I$, see for instances \cite{Cha15,HT18}. The other is to bound or provide explicit formulas for $\Reg~R/I^s$ of some special ideal, see e.g. \cite{Lu21,BHT15}.

The edge ideal  was defined and studied by Villarreal in \cite{Vil90}. Let $G$ be a finite simple graph with  vertex set $V=\{x_1, \cdots, x_n\}$ and  edge set $E$. The ideal generated by all quadratic monomials $x_ix_j$ such that $\{x_i,x_j\}\in E$ is called the {\it edge ideal} of $G$, denoted by $I(G)$. The regularity of powers of edge ideals has been investigated extensively  and  many exciting results appear, see for instance \cite{AB17,BBH19,BBH20,BHT15,JNS18} and the references therein.

 In 1999, Conca and De Negri introduced the concept of a $t$-path ideal in their work \cite{CD99}, which serves as a generalization of the definition of an edge ideal.
 For an integer $2\leq t\leq n$, a {\it  $t$-path} of $G$ is a sequence of $t-1$ distinct edges $\{x_{i_1}, x_{i_2}\}, \{x_{i_2}, x_{i_3}\},\cdots, \{x_{i_{t-1}}, x_{i_t}\}$ of $G$ such that the vertices $x_{i_1}, x_{i_2}, \ldots, x_{i_{t}}$ are pairwise distinct. Such a path is also denoted by $\{x_{i_1}, \ldots, x_{i_t}\}$  for short. The {\it  $t$-path ideal} $I_t(G)$ associated to $G$ is defined as  the squarefree monomial ideal $$I_t(G)=(x_{i_1}\cdots x_{i_t}\mid\{x_{i_1}, \ldots, x_{i_t}\} \textrm{~~is~~ a}~~t\textrm{-path~~ of}~~ G)$$ in the polynomial ring $R$. In the recent years, some algebraic properties of path ideals have been investigated extensively, see for instance \cite{AF18,BHO11,BHO12,HV10,KO14,KS14}.

However, very  little is known about the powers of $t$-path ideals for $t\geq 3$.  In this note, we will specifically concentrate on the powers of path ideals of line graphs.  A line graph of length $n-1$, denote by $L_n$, is a graph with a vertex set $\{x_1, \ldots, x_n\}$ and  an edge set  $\{\{x_j,x_{j+1}\}\mid j=1,\ldots, n-1\}$. That is, $L_n$ has the following form:

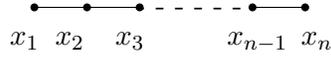
\begin{figure}[ht!]
\begin{tikzpicture}
\draw (4.7,5)--(5.4,5)--(6.1,5);
\draw (7.6,5)--(8.3,5);
\draw (4.9,4.8) node[anchor=north east]{$x_1$};
\draw (5.5,4.8) node[anchor=north east]{$x_2$};
\draw (6.3,4.8) node[anchor=north east]{$x_3$};
\draw (7.6,5.2) node[anchor=north east]{- - - - - -};
\draw (8.2,4.8) node[anchor=north east]{$x_{n-1}$};
\draw (8.8,4.8) node[anchor=north east]{$x_n$};
\fill [color=black] (4.7,5) circle (1.5pt);
\fill [color=black] (5.4,5) circle (1.5pt);
\fill [color=black] (6.1,5) circle (1.5pt);
\fill [color=black] (7.6,5) circle (1.5pt);
\fill [color=black] (8.3,5) circle (1.5pt);
\end{tikzpicture}
\end{figure}
\noindent Recently, ${\rm B\breve{a}l\breve{a}nescu}$-${\rm Cimpoea\textrm{\c{s}}}$ \cite{BC23}  derived an explicit formula for calculating the depth of powers of path ideals associated with line graphs. In contrast, this paper delves into the conditions that govern whether all or some powers of the $t$-path ideals of a line graph exhibit a linear resolution, as detailed in Theorem \ref{main1.1}. Additionally, we establish a regularity formula for powers of the same category of ideals, as presented in Theorem \ref{main1}.

Throughout this note, we denote $L_n$ as the line graph of length $n-1$ and $I_t(L_n)$ as the $t$-path ideal of $L_n$. Let $G(I_t(L_n))$ denote the minimal monomial generating set of $I_t(L_n)$. Namely, $G(I_t(L_n))=\{u_1,\ldots,u_{n-t+1}\}$, where $u_i=x_ix_{i+1}\cdots x_{i+t-1}$ for $i=1, \ldots, n-t+1$.
\smallskip

\section{When powers of  $I_t(L_n)$ have a linear resolution?}

In this section, we study when the powers of  $I_t(L_n)$ have a linear resolution and compute  graded Betti numbers  for such ideals.

 Let $I$ be a monomial ideal that is generated in a single degree and denote by $G(I)$ the minimal set of monomial generators of $I$. Recall from \cite{L} that $I$ is said to be {\it  quasi-linear} if for each $u\in G(I)$, the colon ideal $I\setminus_u:u$ is generated by variables. Here, $I\setminus_u$ is the ideal generated by monomials in $G(I)\setminus \{u\}.$  It was shown in \cite{L} that if $I$ has a linear resolution then $I$ is quasi-linear. Hence we have the following implications:
\begin{equation*}
\begin{aligned}I \mbox{ has linear quotients} \Longrightarrow I \mbox{ has a linear resolution for all field } \K \\ \Longrightarrow I \mbox{ has a linear resolution for some field } \K \Longrightarrow I\mbox{ is quasi-linear}.\end{aligned}
\end{equation*}
When focusing on powers of $t$-path ideals of line graphs, we can demonstrate that the converse of the aforementioned implications also hold true.

We need the following formula that is well-known in the combinatorial theory:
\begin{equation}\label{C}
|\{(a_1,\cdots,a_k)\in \mathbb{Z}_{\geq 0}^k\mid a_1+\cdots+a_k=s \}|=\begin{pmatrix} s+k-1\\k-1\end{pmatrix}.
\end{equation}
Since we can not find a reference for this formula, we give a brief proof utilizing the ``Stars and Bars" method. We envision this problem as placing $k-1$ bars among $s$ stars. For instance, a possible arrangement for $s=7$ and $k=3$ could be \quad $**|*|****$, which corresponds to $a_1=2, a_2=1$ and $a_3=4$. To count the number of such arrangements, we consider choosing $k-1$ positions out of $s+k-1$ positions (comprising of $s$ stars and $k-1$ bars) to place the bars. This can be done using the binomial coefficient $\begin{pmatrix} s+k-1\\k-1\end{pmatrix}$.

\begin{lem} \label{3.1}  Let $\alpha=u_1^{a_1}u_2^{a_2}\cdots u_{n-t+1}^{a_{n-t+1}}$ and $\beta=u_1^{b_1}u_2^{b_2}\cdots u_{n-t+1}^{b_{n-t+1}}$ be monomials in $G(I_t(L_n)^s)$. Then $\alpha=\beta$ if and only if $a_i=b_i$ for $i=1,\ldots,n-t+1$. Consequently, we have $$|G(I_t(L_n)^s)|=\begin{pmatrix} s+n-t\\s\end{pmatrix}.$$
\end{lem}

\begin{pf} To prove the first assertion, it suffices to prove if  $\alpha=\beta$ then $a_i=b_i$ for all $i$. We proceed by induction on $s$. The base case when $s=1$ is straightforward. Now, assume $s>1$. Let $k$ be the smallest  $i$  such that $x_i$ divides $\alpha$ (since $\alpha=\beta$, it also divides $\beta$). Then $a_k>0$, $b_k>0$ and $a_i=b_i=0$ for $i=1,\cdots,k-1$. From this, it follows that $$\frac{\alpha}{u_k}=u_k^{a_k-1}u_{k+1}^{a_{k+1}}\cdots u_{n-t+1}^{a_{n-t+1}}=\frac{\beta}{u_k}=u_k^{b_k-1}u_{k+1}^{b_{k+1}}\cdots u_{n-t+1}^{b_{n-t+1}}.$$ Hence,  by induction hypothesis, we have $a_k-1=b_k-1$ and $a_i=b_i$ for $i=k+1,\ldots,n-t+1$.
This completes the proof of the first assertion.

  The second statement follows directly from Formula~\ref{C}.
\end{pf}

\vspace{2mm}
To facilitate the proof of the subsequent result, we introduce additional notation. Assume  that $\alpha,\beta$ are monomials in $R$. We define $$\deg_k\alpha:=\Max\{i\mid x_k^i\mbox{ divides }\alpha\} \mbox{ for } k=1,\ldots, n, \mbox{ and } \deg\alpha:=\sum_{k=1}^{n}\deg_k\alpha.$$  In addition, we use $\alpha/\beta$ to denote the monomial $\frac{\alpha}{(\alpha,\beta)}$. Here, $(\alpha,\beta)$ is the greatest common divisor of $\alpha$ and $\beta$.  We now present the first main result of this paper.

\begin{theorem}\label{main1.1}
Let $I=I_t(L_n)$ be the $t$-path ideal of $L_n$. Then the following statements are equivalent:
\begin{enumerate}
\item $I^s$ has a linear resolution for some $s\geq 1$;
\item $I^s$ has a linear resolution for all $s\geq 1$;
\item $t\leq n\leq 2t$;
\item $I^s$ has linear quotients for some $s\geq 1$;
\item $I^s$ has linear quotients for all $s\geq 1$;
\item $I^s$ is quasi-linear for some $s\geq 1$;
\item $I^s$ is quasi-linear for all $s\geq 1$.
\end{enumerate}
\end{theorem}

\begin{pf} The implications  $(2)\Rightarrow (1)$, $(5)\Rightarrow (4)$ and $(7)\Rightarrow (6)$ are automatical.

The implications $(4)\Rightarrow (1)$ and $(5)\Rightarrow (2)$  are derived from \cite[Proposition 8.2.1]{HH11}.

The implications $(1)\Rightarrow (6)$ and $(2)\Rightarrow (7)$   are consequences of  \cite[Theorem 2.3]{L}.

Therefore,  we only need to prove $(6)\Rightarrow (3)\Rightarrow (5)$.

$(6)\Rightarrow (3)$ Assume on the contrary that $n\geq 2t+1$. Fix $s\geq 1$.  Our goal is to  prove $I^s$ is not quasi-linear. To this end, we denote $\alpha=u_{n-t+1}^s$ and let $J$ be the ideal generated by monomials in $G(I^s)\setminus \{\alpha\}$. Choose $k\in [n]$ such that the variable $x_k$ belongs to $J:\alpha$. Then, there is $\beta\in G(I^s)\setminus \{\alpha\}$ such that $x_k=\beta/ \alpha=\frac{\beta}{(\alpha,\beta)}$. Since $\deg \beta=st$, it follows that  $(\alpha,\beta)$ has a degree of $ts-1$. Note that $\alpha=u_{n-t+1}^s=x_{n-t+1}^s\ldots x_{n}^s$.  Consequently, $\beta$ takes the form $x_kx_{n-t+1}^{a_1}\ldots x_{n}^{a_t}$, where exactly one of the $a_i$'s is equal to $s-1$, while all the others are equal to $s$. However, if $a_t=s$, $\beta$ would be equal to $\alpha$, which is a contradiction. This implies $a_t=s-1$ and $\beta=u_iu_{n-t+1}^{s-1}$ for some $1\leq i\leq n-t$. Since $x_{n-t+1}x_{n-t+2}\cdots x_{n-1}$ divides $u_i$, we can deduce that $i=n-t$, and consequently $\beta=u_{n-t}u_{n-t+1}^{s-1}$. Hence, $x_k=x_{n-t}$ and, as a result, $x_{n-t}$ is the unique variable in $J:\alpha$. Noting that $x_{n-t}$ does not divides $u_1^s$, it follows that $J:\alpha$ is not generated solely by variables. Therefore, we have established that $I^s$ is not quasi-linear, as desired.

$(3)\Rightarrow (5)$  By an easy calculation, we have $(u_1,\ldots,u_{i-1}):u_i$ is generated by the unique variable $x_{i-1}$ for $i=2,\ldots,n-t+1$. Hence $I^s$ has linear quotients for $s=1$. As for the situation where $s\geq  2$, we first draw attention to the following observation stemming from the condition $n\leq 2t$:

   \textbf{Fact ($\diamond$}): If $1\leq k\leq n-t$, then
  $\deg_{k} u_1=\deg_{k} u_2=\cdots=\deg_k{u_k}=1$ and $\deg_{k} u_{k+1} =\cdots=\deg_k{u_{n-t+1}}=0$.

 In view of Lemma \ref{3.1}, we have $$G(I^s)=\{ u_1^{a_1}\cdots u_{n-t+1}^{a_{n-t+1}}\mid a_1+\cdots+a_{n-t+1}=s, a_i\in \mathbb{Z}_{\geq 0} \mbox{ for all } i=1, \ldots, n-t+1\}.$$   We define a linear order on $G(I^s)$ as follows:
 if $\alpha=u_1^{a_1}\cdots u_{n-t+1}^{a_{n-t+1}}$ and $\beta=u_1^{b_1}\cdots u_{n-t+1}^{b_{n-t+1}}$ are elements of $G(I^s)$, then $$\alpha<\beta \Longleftrightarrow \mbox{ if there exists } i>0  \mbox{ such that } a_1=b_1,\ldots, a_{i-1}=b_{i-1}\mbox{ and }a_i<b_i.$$ With respect to this order, $u_{n-t+1}^s$ and $u_1^s$ are the least element and greatest element of $G(I^s)$ respectively.

  Denote by $q$ the cardinality of  $G(I^s)$ and let $\alpha_1,\alpha_2,\ldots,\alpha_q$ be all the elements of $G(I^s)$ such that  $\alpha_q>\alpha_{q-1}>\cdots >\alpha_1$. We want to show that for each $1\leq j\leq q-1$, the colon ideal $C_j:=(\alpha_q, \ldots, \alpha_{j+1}): \alpha_j$ is generated by variables.

Fix $1\leq j<q$ and assume $\alpha_j= u_1^{a_1}\cdots u_{n-t+1}^{a_{n-t+1}}$. We may write  $\{1\leq \lambda\leq n-t+1\mid a_\lambda>0\}$ as $\{\lambda_1,\lambda_2,\ldots,\lambda_g\}$ with $\lambda_1<\lambda_2<\cdots<\lambda_g$. We claim that $C_j=(x_{\lambda_1-1},\ldots,x_{\lambda_g-1})$. Here, we use the convention that $x_0=0$. Namely, if $\lambda_1=1$, then $C_j=(x_{\lambda_2-1},\ldots,x_{\lambda_g-1})$.

 Take any $\lambda\in \{\lambda_1,\lambda_2,\ldots,\lambda_g\}$ with $\lambda\geq 2$. Since $\alpha:=u_1^{a_1}\cdots u_{\lambda-1}^{a_{\lambda-1}+1}u_\lambda^{a_\lambda-1}\cdots u_{n-t+1}^{a_{n-t+1}}>\alpha_j$, we have $\alpha/\alpha_j=x_{\lambda-1}\in C_j$. This proves the inclusion  $(x_{\lambda_1-1},\ldots,x_{\lambda_g-1})\subseteq C_j$. Conversely, let $\beta=u_1^{b_1}\cdots u_{n-t+1}^{b_{n-t+1}}$ with $\beta>\alpha_j$. Then there exists $1\leq \ell< \lambda_g$ such that $a_1=b_1$, $\ldots$, $a_{\ell-1}=b_{\ell-1}$ and $a_{\ell}<b_{\ell}$.  If $\ell<\lambda_1$,  we have $\deg_{(\lambda_1-1)}\beta\geq b_{\ell}>\deg_{(\lambda_1-1)}\alpha_j=0$ due to \textbf{Fact ($\diamond$)}. Hence $\beta/\alpha_j\in(x_{\lambda_1-1})\subseteq (x_{\lambda_1-1},\ldots,x_{\lambda_g-1})$. Suppose next that $\lambda_f\leq \ell <\lambda_{f+1}$ for some $f\in\{1,2, \ldots, g-1\}$. Then, taking into account \textbf{Fact ($\diamond$)},  we have $$\deg_{(\lambda_{f+1}-1)}\beta=b_1+\cdots+b_{(\lambda_{f+1}-1)}> a_1+\cdots+a_{(\lambda_{f+1}-1)}=\deg_{(\lambda_{f+1}-1)}\alpha_i,$$ and so $x_{(\lambda_{f+1}-1)}$ divides $\beta/\alpha_j$. It follows that  $\beta/\alpha_j\in(x_{(\lambda_{f+1}-1)})\subseteq (x_{\lambda_1-1},\ldots,x_{\lambda_g-1})$ and thus the claim follows. This completes the proof. \end{pf}

  It is a known fact that if an ideal $I$ possesses  a linear resolution, then $\Reg~I$ is equal to the degree of minimal generators of $I$. Consequently, the following result is immediate.

\begin{coro} \label{3.2}
Let $I=I_t(L_n)$ be the $t$-path ideal of $L_n$. If $t\leq n\leq 2t$, then $\Reg~R/I^s=ts -1$.
\end{coro}

We now compute the Betti numbers of $I_t(L_n)^s$ for  $t\leq n\leq 2t$ by employing \cite[Corollary 8.2.2.]{HH11}. According to this corollary, if $I$ is generated by a sequence $f_1,\ldots,f_p$ of monomials all having the same degree, and the colon ideal $(f_1,\ldots,f_{j-1}):f_j$ is generated by $r_j$ variables for each $j=2,\ldots,p$, then the Betti number $\beta_i(I)$ can be determined as $\sum_{j=2}^p\binom{r_j}{i}$.
\begin{coro}
Let $I=I_t(L_n)$ be the $t$-path ideal of $L_n$. If $t\leq n\leq 2t$, then $$\beta_i(I^s)=\beta_{i,i+st}(I^s)=\sum_{k=i}^{n-t}\begin{pmatrix}n-t\\k\end{pmatrix} \begin{pmatrix}s\\k\end{pmatrix}\begin{pmatrix}k\\i\end{pmatrix}.$$

In particular, $\Pd~R/I^s=\min\{n-t+1, s+1\}$.
\end{coro}
\begin{pf}
For $1\leq j<q$, let $C_j$ represent the same ideal as defined in the proof of $(3)\Rightarrow (5)$ in Theorem~\ref{main1.1}, and let $r_j$ be the number of variables in $C_j$. According to the proof of $(3)\Rightarrow (5)$ in Theorem~\ref{main1.1}, if $\alpha_j = u_1^{a_1}\cdots u_{n-t+1}^{a_{n-t+1}}$, then $C_j$ is generated by the variables $x_{\lambda-1}$ for which $a_{\lambda} > 0$ and $\lambda \geq 2$. Consequently, $r_j$ is equal to the number of $\{2\leq \lambda\leq n-t+1\mid a_{\lambda} > 0\}$. It is evident that $r_j$ falls within the range $1\leq r_j\leq n-t$.

 For each $1\leq k\leq n-t$, let $S_k$ denote  the number of $\{1\leq j< q\mid r_j=k \}$. Note that $\alpha_q=u_1^{s}$. It follows that $S_k$   is equal to  the number of the set $A_k$, which, by definition,  is composed of all the tuples $(a_1, a_2,\ldots, a_{n-t+1})$  such that  $a_1+\cdots+a_{n-t+1}=s$, $a_1\neq s$ and $|\{2\leq \lambda\leq n-t+1\mid a_{\lambda}>0\}|=k.$ Since $k\geq 1$, condition that $a_1\neq s$ can be skipped. To compute the $S_k$, we first compute the number of the set $A_{k,e}:=\{(a_1,\ldots,a_{n-t+1})\in A_k\mid a_1=e\}$ for $e=0,\ldots, s-k$. It is not difficult to see that $|A_{k,e}|$ is equal to  $\binom{n-t}{k}$ multiplied by the number of the positive integral solutions of the equation $b_1+\cdots+b_k=s-e$. Furthermore, the number is equal to $\binom{s-e-1}{k-1}$ by an application of  Formula~\ref{C}.   This implies $$S_k=\sum_{e=0}^{s-k}|A_{k,e}|=\begin{pmatrix}n-t\\k\end{pmatrix}\left (\begin{pmatrix} s-1\\ k-1\end{pmatrix}+\cdots+ \begin{pmatrix} k-1\\ k-1\end{pmatrix}\right )=\begin{pmatrix}n-t\\k\end{pmatrix}\begin{pmatrix}s\\k\end{pmatrix}.$$
It follows  from \cite[Corollary 8.2.2.]{HH11} that $$\beta_i(I^s)=\sum_{j=1}^{q-1}\binom{r_j}{i}= S_i\begin{pmatrix}i\\i\end{pmatrix}+ S_{i+1}\begin{pmatrix}i+1\\i\end{pmatrix}+ \cdots+S_{n-t}\begin{pmatrix}n-t\\i\end{pmatrix}.$$
The last statement  follows since $\Pd~R/I^s=\Pd~I^s+1$ and $\Pd~I^s=\max\{i\mid \beta_{i}(I^s)\neq 0\}$.
\end{pf}

\section{Regularity of powers of $I_t(L_n)$}

\setcounter{equation}{0}
\renewcommand{\theequation}
{3.\arabic{equation}}

In this section we give an explicit formula for the regularity of the powers of $I_t(L_n)$.

We first recall a result in \cite{AF18} by  Alilooee and Faridi.  Before doing so we define a function $\Gamma$ on $\mathbb{Z}_{\geq 0}\times \mathbb{Z}_{\geq 0}$ as follows: $$\Gamma(n,t):=\left\{
        \begin{array}{rcl}
         p(t-1),\;\;\;\;\;\;\;\;     & & \mbox{if } n=p(t+1)+d \mbox{  and  }  0\leq d< t;  \\
         \\
         (p+1)(t-1), & & \mbox{if  } n=p(t+1)+t.
         \end{array}
         \right. .$$

         In light of this notation, the second part of \cite[Corollary 4.15]{AF18} can be restated as follows.
         \begin{lem}\label{3.0} Let $n,t,p,d$ be postive integers such that $n=p(t+1)+d$ with $0\leq d\leq t$. Then $\Reg~R/I_t(L_n)=\Gamma(n,t).$
         \end{lem}

It is worth noting that \cite[Corollary 4.15]{AF18} states that $\Reg~R/I_t(L_n)=\Gamma(n,t)$ when $n\geq t\geq 2$. However, this equality actually holds true for  all $n\geq 1$ and $t\geq 1$, as stated in Lemma~\ref{3.0}. In fact, when either $n<t$ or $t=1$ holds, we observe that $\Reg~R/I_t(L_n)=\Gamma(n,t)=0$. The following lemma collects some easy but useful observations on  the function $\Gamma(n,t)$.
\begin{lem}\label{f}
Let $n,t,a,b$ be  integers $\geq 0$ with $n\geq t+1$. Then we have\\
{\rm (1)} $\Gamma(n-t-1,t)=\Gamma(n,t)-(t-1)$;\\
{\rm (2)} $\Gamma(a,t)+\Gamma(b,t)\leq \Gamma(a+b+1,t)$ for all $a,b\geq 1$.
\end{lem}
\begin{pf}
{\rm (1)} Straightforward.

{\rm (2)} Let $a=p_1(t+1)+d_1$ and $b=p_2(t+1)+d_2$, where $p_1,p_2\geq 0$ and $0\leq d_1,d_2\leq t$. As a consequence we have $a+b+1=(p_1+p_2)(t+1)+(d_1+d_2+1)$.

If both $d_1$ and $d_2$ are equal to $t$. Then $\Gamma(a,t)=(p_1+1)(t-1)$, $\Gamma(b,t)=(p_2+1)(t-1)$. Note that $a+b+1=(p_1+p_2+1)(t+1)+t$. Then $\Gamma(a+b+1,t)=(p_1+p_2+2)(t-1)=\Gamma(a,t)+\Gamma(b,t)$.

If exactly one of  $d_1$ and $d_2$ is equal to $t$. We may assume that $d_1=t$ and $d_2\neq t$. Then $\Gamma(a,t)=(p_1+1)(t-1)$, $\Gamma(b,t)=p_2(t-1)$. Note that $a+b+1=(p_1+p_2+1)(t+1)+d_2$. Then $\Gamma(a+b+1,t)=(p_1+p_2+1)(t-1)=\Gamma(a,t)+\Gamma(b,t)$.

If neither $d_1$ nor $d_2$ is equal to $t$, then $\Gamma(a,t)=p_1(t-1)$ and $\Gamma(b,t)=p_2(t-1)$. Note that $$a+b+1=\left\{
        \begin{array}{rcl}
         (p_1+p_2)(t+1)+(d_1+d_2+1),\;\;\;\;\;    & &  d_1+d_2+1\leq t  \\
         \\
         (p_1+p_2+1)(t+1)+(d_1+d_2-t),     & &  d_1+d_2+1\geq t+1
         \end{array}
         \right.,$$
it follows that  $$\Gamma(a+b+1,t)=\left\{
        \begin{array}{rcl}
         (p_1+p_2)(t-1),\;\;\;\;\;\;& &  d_1+d_2+1< t  \\
         \\
         (p_1+p_2+1)(t-1),          & &  d_1+d_2+1\geq t
         \end{array}
         \right. .$$
Hence $\Gamma(a,t)+\Gamma(b,t)\leq \Gamma(a+b+1,t)$. This finishes the proof.
\end{pf}

We record the following known results for the later use.
\begin{lem} {\cite[Lemma 3.1]{HT100}} \label{reg}
Let $0\ra M\ra N\ra P\ra 0$ be a short exact sequence of finitely generated graded $S$-modules. Then\\
{\rm (1)} $\Reg~ N\leq{\rm max}\{\Reg~ M, \Reg~ P\}$. The equality holds for $\Reg~P\neq \Reg~M-1$.\\
{\rm (2)} $\Reg~ M\leq{\rm max}\{\Reg~ N, \Reg~ P+1\}$. The equality holds for $\Reg~N\neq \Reg~P$.\\
{\rm (3)} $\Reg~ P\leq{\rm max}\{\Reg~ M-1, \Reg~ N\}$. The equality holds for $\Reg~M\neq \Reg~N$.
\end{lem}

\begin{lem}\label{splitting}{\cite[Lemma 3.2]{HT100}}
Let $A=\K[\textbf{x}]=\K[x_1,\ldots, x_n]$, $B=\K[\textbf{y}]=\K[y_1,\ldots,y_N]$ and $S=\K[\textbf{x}, \textbf{y}]$ be polynomial rings. Then for nonzero homogeneous ideals $I\subset A$ and $J\subset B$, we have \\
{\rm (1)} $\Reg~S/(I+J)=\Reg~A/I+\Reg~B/J$;\\
{\rm (2)} $\Reg~S/(IJ)=\Reg~A/I+\Reg~B/J+1$.
\end{lem}

Let $R$ be a polynomial ring over a field. For a monomial $u\in R$, the {\it support} of $u$, denoted by $\mbox{supp}(u)$, is the set of variables of $R$ that divide $u$.
Let $I$ be a monomial ideal of $R$. Then the {\it support}  of $I$, denoted as $\mbox{supp}(I)$, is the union of supports of minimal monomial generators of $I$. According to Lemma~\ref{splitting}, if $I$ and $J$ are monomial ideals of $R$ such that $\mbox{supp}(I)\cap \mbox{supp}(J)=\emptyset$, then $\Reg~R/(I+J)=\Reg~R/I+\Reg~R/J$.

\begin{lem}{\cite[Lemma 2.1]{BC23}}\label{s-1}
Let$t,n$ and $s$ be integers such that $1\leq t\leq n$ and $s\geq 2$. Then $$I_t(L_n)^s: u_{n-t+1}=I_t(L_n)^{s-1}.$$
\end{lem}

Let $I,J$ and $K$ be monomial ideals. It is known that $(I+J):K=I:K+J:K$. We shall utilize this result in our subsequent discussion without explicit reference.
We now are ready to prove the second main result of this paper.
\begin{theorem}\label{main1}
 Assume $n\geq t\geq 2$. Let $I_t(L_n)$ be the $t$-path ideal of $L_n$. Then, for all $s\geq 1$, the regularity of $I_t(L_n)^s$ is given by $$\Reg~R/I_t(L_n)^s=\Gamma(n,t)+t(s-1).$$
\end{theorem}
\begin{pf} We proceed by induction on $s$.
If  $s=1$,  the result is given by Lemma~\ref{3.0} (i.e., \cite[Corollary 4.15]{AF18}). Suppose now $s\geq 2$. Since the result is clear from Corollary~\ref{3.2} in the case that $t\leq n\leq 2t$, we may assume that $n\geq 2t+1$ from now on.

 As before we write $$I_t(L_n)=(u_1, u_2,\ldots ,u_{n-t+1}),$$ where $u_i$ denotes the monomial $x_ix_{i+1}\cdots x_{i+t-1}$ for $1\leq i\leq n-t+1$. We consider the following short exact sequence: $$0 \ra\frac{R}{(I_t(L_n)^s:u_{n-t+1})}[-t] \mathop{\ra}\limits^{\cdot u_{n-t+1}}\frac{R}{I_t(L_n)^s} \ra\frac{R}{(I_t(L_n)^s,u_{n-t+1})} \ra 0.$$ By induction hypothesis and by Lemma \ref{s-1},  we have $$\Reg~\frac{R}{(I_t(L_n)^s: u_{n-t+1})}=\Reg~\frac{R}{I_t(L_n)^{s-1}}=\Gamma(n,t)+t(s-2).$$ Hence $$\Reg~\frac{R}{(I_t(L_n)^s: u_{n-t+1})}[-t]=\Reg~\frac{R}{(I_t(L_n)^s: u_{n-t+1})}+t=\Gamma(n,t)+t(s-1).$$
      We claim that $\Reg~\frac{R}{(I_t(L_n)^s,u_{n-t+1})}= \Gamma(n,t)+t(s-1)$. To prove this claim, we put $$A_0=(I_t(L_n)^s,u_{n-t+1}),$$ and   then, for $1\leq j\leq n-t$, we put $$A_j=(A_{j-1}, u_{n-t-j+1}) \mbox{ and } B_j=A_{j-1}: u_{n-t-j+1}$$  recursively.  It follows that
\begin{equation}\label{A_j}
A_j=(A_0, u_{n-t},\ldots, u_{n-t-j+1})=(I_t(L_n)^s, u_{n-t+1}, \ldots, u_{n-t-j+1}).
\end{equation}
In particular, we have  $A_{n-t}=(I_t(L_n)^s, u_{n-t+1},\ldots, u_1)=I_t(L_n)$, and thus
\begin{equation}\label{A_{n-t}}
\Reg~\frac{R}{A_{n-t}}=\Gamma(n,t).
\end{equation}
According to Lemma~\ref{s-1}, for any $1\leq j\leq n-t$,  we have
\begin{eqnarray}
  B_j &=& (I_t(L_n)^s, u_{n-t+1},\ldots, u_{n-t-j+2}): u_{n-t-j+1}\nonumber \\
      &=& ((u_1,\ldots,u_{n-t-j+1})^s, u_{n-t+1}, \ldots, u_{n-t-j+2}): u_{n-t-j+1} \nonumber \\
      &=& ((u_1,\ldots,u_{n-t-j+1})^{s-1}, u_{n-t+1}, \ldots, u_{n-j+2}, x_{n-j+1}) \nonumber \\
      &=& I_t(L_{n-j})^{s-1}+(u_{n-t+1}, \ldots, u_{n-j+2})+(x_{n-j+1}).
\end{eqnarray}According to the induction hypothesis, we have $\Reg~\frac{R}{I_t(L_{n-j})^{s-1}}=\Gamma(n-j,t)+t(s-2)$. On the other hand, it is noteworthy that the ideal $(u_{n-t+1}, \ldots, u_{n-j+2})$ is isomorphic to the $t$-path ideal of a line graph of length $j-1$.   Since  the monomial ideals $I_t(L_{n-j})^{s-1}$, $(u_{n-t+1}, \ldots, u_{n-j+2})$ and $(x_{n-j+1})$ have  pairwise disjoint supports,
 applying Lemma \ref{splitting}(1) and Lemma \ref{f}(2), we can conclude that
\begin{eqnarray}\label{B_j}
  \Reg~\frac{R}{B_j}[-t] &=& \Reg~\frac{R}{I_t(L_{n-j})^{s-1}}+\Reg~\frac{R}{I_t(L_{j-1})}+t \nonumber \\
   &=& \Gamma(n-j,t)+t(s-2)+\Gamma(j-1,t)+t \nonumber \\
   &\leq& \Gamma(n,t)+t(s-1).
\end{eqnarray}
We now show that (\ref{B_j}) is actually  an equality when $j=n-t$.
In fact, we have $$B_{n-t}=I_t(L_t)^{s-1}+I_t(L_{n-t-1})+(x_{t+1}).$$  From this it follows that
\begin{eqnarray}\label{B_{n-t}}
  \Reg~\frac{R}{B_{n-t}}[-t] &=& \Reg~\frac{R}{I_t(L_t)^{s-1}}+\Reg~\frac{R}{I_t(L_{n-t-1})}+t \nonumber\\
   &=& t(s-1)-1+\Gamma(n-t-1,t)+t \nonumber\\
   &=& t(s-1)-1+\Gamma(n,t)-(t-1)+t \nonumber\\
   &=& \Gamma(n,t)+t(s-1).
\end{eqnarray}
Here, the third equality follows from Lemma~\ref{f}(1). Now, we consider the following short exact sequences of  cyclic graded  $R$-modules: $$0 \ra\frac{R}{B_1}[-t] \mathop{\ra}\limits^{\cdot u_{n-t}}\frac{R}{A_0} \ra\frac{R}{A_1} \ra 0;$$
$$0 \ra\frac{R}{B_2}[-t] \mathop{\ra}\limits^{\cdot u_{n-t-1}}\frac{R}{A_1} \ra\frac{R}{A_2} \ra 0;$$
$$\vdots$$
$$0 \ra\frac{R}{B_{n-t}}[-t] \mathop{\ra}\limits^{\cdot u_1}\frac{R}{A_{n-t-1}} \ra\frac{R}{A_{n-t}} \ra 0.$$
 Considering (\ref{A_{n-t}})  together with (\ref{B_{n-t}}),  and applying Lemma~\ref{reg}(1) to  the last  short exact sequence mentioned above,  we obtain \begin{eqnarray}\label{3.6} \Reg~\frac{R}{A_{n-t-1}}=\Gamma(n,t)+t(s-1).  \end{eqnarray}
Considering (\ref{3.6}) together with (\ref{B_j}), we have $\Reg \frac{R}{B_{n-t-1}}[-t]\leq \Gamma(s,t)+t(s-1)=\Reg \frac{R}{A_{n-t-1}}$. In particular, $\Reg \frac{R}{A_{n-t-1}}\neq \Reg \frac{R}{B_{n-t-1}}[-t]-1$. Taking into account these facts  and by applying Lemma \ref{reg}(1) to the the last second short exact sequence above, we obtain
 $$\Reg~\frac{R}{A_{n-t-2}}=\Gamma(n,t)+t(s-1).$$
 Continuing in this way, finally  we arrive at  $$\Reg~\frac{R}{A_0}=\Gamma(n,t)+t(s-1).$$ This proves the claim and so the proof is complete.
\end{pf}

\begin{remark}\em
By checking the proof of Theorem \ref{main1}, we  obtain  not only the regularity of $\frac{R}{I_t(L_n)^s}$,  but also  the following  formula:
   $$\Reg~\frac{R}{(I_t(L_n)^s, u_{n-t+1}, \ldots, u_j)}=\Gamma(n,t)+t(s-1),$$
for   all $s\geq 1$ and $2\leq j\leq n-t+1$.
\end{remark}

{\bf \noindent Acknowledgment:}
This research is supported by NSFC (No. 11971338).  We would express our sincere gratitude to the referee for his/her careful reading and numberous comments that improve the presentation of our paper greatly.

\bibliography{}

\begin{thebibliography}{9999}
\bibitem{AF18} A. Alilooee and S. Faridi, \textit{Graded Betti numbers of path ideals of cycles and lines}, J. Algebra Appl. \textbf{17}(2018), 1850011, 1--17.  \par
\bibitem{AB17} A. Alilooee and A. Banerjee, \textit{Powers of edge ideals of regularity three bipartite graphs}, J. Commut. Algebra \textbf{9}(2017),441--454.  \par
\bibitem{BBH19} A. Banerjee, S. K. Beyarslan and H. T. ${\rm H\grave{a}}$, \textit{Regularity of edge ideals and their powers}, Advances in Algebra, 17--52, Springer Proc. Math. Stat. \textbf{277}, Springer, Cham, 2019.\par
\bibitem{BBH20} A. Banerjee, S. K. Beyarslan, and H. T. ${\rm H\grave{a}}$, \textit{Regularity of powers of edge ideals: from local properties to global bounds}. Algebr. Comb. \textbf{3}(2020),839--854.\par
\bibitem{BC23} S. ${\rm B\breve{a}l\breve{a}nescu}$, M. $\textrm{Cimpoea\c{s}}$, \textit{Depth and Stanley depth of powers of the path ideal of a path graph}, preprint, 2023.\par
\bibitem{BHO11} R. R. Bouchat, H. T. ${\rm H\grave{a}}$ and A. O'Keefe, \textit{Path ideals of rooted trees and their graded Betti numbers}, J. Combin. Theory Ser. A \textbf{118}(2011) 2411--2425.\par
\bibitem{BHO12} R. R. Bouchat, H. T. ${\rm H\grave{a}}$ and A. O'Keefe, \textit{Corrigendum to "Path ideals of rooted trees and their graded Betti numbers" [J. Combin. Theory Ser. A \textbf{118}(2011) 2411--2425]}, J. Combin. Theory Ser. A \textbf{119}(2012) 1610--1611.\par
\bibitem{BHT15} S. K. Beyarslan, H. T. ${\rm H\grave{a}}$ and T. N. Trung, \textit{Regularity of powers of forests and cycles}, J. Algebraic Combin. \textbf{42}(2015), 1077--1095.\par
\bibitem{Cha15} M. Chardin, \textit{Regularity stabilization for the powers of graded M-primary ideals}, Proc. Amer. Math. Soc. \textbf{143}(2015), 3343--3349.\par
\bibitem{CD99} A. Conca, E. De Negri, \textit{M-sequences, graph ideals and ladder ideals of linear types}, J. Algebra \textbf{211}(1999), 599--624.\par
\bibitem{CHT99} S. D. Cutkosky, J. Herzog, N. V. Trung, \textit{Asymptotic Behaviour of the Castelnuovo-Mumford Regularity}, Compositio Math. \textbf{118}(1999), 243--261.\par
\bibitem{HH11} J. Herzog, T. Hibi, \textit{Monomial Ideals}, Graduate Texts in Mathematics, vol. 260. Springer, Berlin, 2011.\par
\bibitem{HT100} L. T. Hoa, N. D. Tam, \textit{On some invariants of a mixed product of ideals}, Arch. Math. \textbf{94}(2010), 327--337.\par
\bibitem{HT18} N. T. Hang, T. N. Trung, \textit{Regularity of powers of cover ideals of unimodular hypergraphs}, J. Algebra \textbf{513}(2018), 159--176.\par
\bibitem{HV10} J. He, A. Van Tuyl, \textit{Algebraic properties of the path ideal of a tree}, Comm. Algebra \textbf{38}(2010), 1725--1742.\par
\bibitem{JNS18} A. V. Jayanthan, N. Narayanan and S. Selvaraja, \textit{Regularity of powers of bipartite graphs}, J. Algebraic Combin. \textbf{47}(2018), 17--38. \par
\bibitem{Kod00} V. Kodiyalam, \textit{Asymptotic behaviour of Castelnuovo-Mumford regularity}, Proc. Amer. Math. Soc. \textbf{128}(2000), 407--411.\par
\bibitem{KO14} M. Kubitzke, A. Olteanu, \textit{Algebraic properties of classes of path ideals of posets}, J. Pure Appl. Algebra \textbf{218}(2014), 1012--1033.\par
\bibitem{KS14} D. Kiani, S. Saeedi Madani, \textit{Betti numbers of path ideals of trees}, Comm. Algebra \textbf{44}(2016) 5376--5394. \par
\bibitem{Lu21} D. Lu, \textit{Geometric regularity of powers of two-dimensional
squarefree monomial ideals}, J. Algebraic Combin. \textbf{53}(2021), 991--1014. \par
\bibitem{L} D. Lu, \textit{Linear resolutions and quasi-linearity of monomial ideals},  Mathematica Scandinavica,  \textbf{129}(2023) 189--208. \par
\bibitem{Vil90} R. H. Villarreal, \textit{Cohen-Macaulay graphs}, Manuscripta Math. \textbf{66}(1990), 277--293.\par
\end{thebibliography}

\end{document}